\documentclass[graybox]{svmult}

\usepackage{mathrsfs}
\usepackage{mathtools}
\usepackage{todonotes}

\usepackage{pgf,tikz}
\usepackage{mathrsfs}
\usetikzlibrary{arrows}

\usepackage{mathptmx}       
\usepackage{helvet}         
\usepackage{courier}        
\usepackage{type1cm}        
\usepackage{makeidx}         
\usepackage{graphicx}        
\usepackage{multicol}        
\usepackage[bottom]{footmisc}

\usepackage{cite}

\usepackage{amsmath,amssymb,bm}

\makeindex             

\begin{document}
\label{chapter:HOIM-I}

\title*{High-Order Isogeometric Methods for Compressible Flows.} 
\subtitle{I. Scalar Conservation Laws}
\titlerunning{High-Order Isogeometric Methods for Compressible Flows. I.}
\author{Andrzej Jaeschke and Matthias M{\"o}ller}
\institute{Andrzej Jaeschke \at \L{}\'od\'z University of Technology, Institute of Turbomachinery, ul. W\'olcza\'nska 219/223, 90-924 \L{}\'od\'z, Poland, \email{andrzej.jaeschke@p.lodz.pl}
\and Matthias M{\"o}ller \at Delft University of Technology, Faculty of Electrical Engineering, Mathematics and Computer Science, Delft Institute of Applied Mathematics, Van Mourik Broekmanweg 6, 2628 XE Delft, The Netherlands, \email{m.moller@tudelft.nl}}
%

%
\maketitle

\abstract*{Isogeometric analysis was applied very successfully to many
  problem classes like linear elasticity, heat transfer and
  incompressible flow problems but its application to compressible
  flows is very rare. However, its ability to accurately represent
  complex geometries used in industrial applications makes IGA a
  suitable tool for the analysis of compressible flow problems that
  require the accurate resolution of boundary layers. The
  convection-diffusion solver presented in this chapter, is an
  indispensable step on the way to developing a compressible flow
  solver for complex viscous industrial flows. It is well known that
  the standard Galerkin finite element method and its isogeometric
  counterpart suffer from spurious oscillatory behaviour in the
  presence of shocks and steep solution gradients. As a remedy, the
  algebraic flux correction paradigm is generalized to B-Spline basis
  functions to suppress the creation of oscillations and occurrence of
  non-physical values in the solution. This work provides early
  results for scalar conservation laws and lays the foundation for
  extending this approach to the compressible Euler equations in
  \cite{igaAFC2_fef2017}.
}

\abstract{Isogeometric analysis was applied very successfully to many
  problem classes like linear elasticity, heat transfer and
  incompressible flow problems but its application to compressible
  flows is very rare. However, its ability to accurately represent
  complex geometries used in industrial applications makes IGA a
  suitable tool for the analysis of compressible flow problems that
  require the accurate resolution of boundary layers. The
  convection-diffusion solver presented in this chapter, is an
  indispensable step on the way to developing a compressible solver
  for complex viscous industrial flows. It is well known that the
  standard Galerkin finite element method and its isogeometric
  counterpart suffer from spurious oscillatory behaviour in the
  presence of shocks and steep solution gradients. As a remedy, the
  algebraic flux correction paradigm is generalized to B-Spline basis
  functions to suppress the creation of oscillations and occurrence of
  non-physical values in the solution. This work provides early
  results for scalar conservation laws and lays the foundation for
  extending this approach to the compressible Euler equations in
  \cite{igaAFC2_fef2017}.
}

\section{Introduction}
\label{sec:introduction}

Isogeometric analysis (IGA) was proposed by Hughes et al. in \cite{IGAArticle}. Since its birth it was successfully applied in a variety of use case scenarios ranging from linear elasticity and incompressible flows to fluid-structure interaction problems \cite{IGABook}. There were, however, not many approaches to apply this method to compressible flow problems \cite{IGAart2, Jaeschke}. Although this application did not gain the attention of many researches yet, it seems to be a promising field. Flow problems are usually defined on domains with complex but smooth shapes, whereby the exact representation of the boundary is indispensable due to the crucial influence of boundary layers on the flow behaviour. This is where IGA has the potential to demonstrate its strengths.

It is a well know fact that standard Galerkin finite element schemes (FEM) suffer from infamous instabilities when applied to convection-dominated problems, such as compressible flows. The same unwanted behaviour occurs for IGA-based standard Galerkin schemes \cite{IGABook} making it necessary to develop high-resolution high-order isogeometric schemes that overcome these limitations. From the many available approaches including the most commonly used ones, i.e., the streamline upwind Petrov–Galerkin (SUPG) method introduced by Brooks and Hughes in \cite{Brooks1982199}, we have chosen for the algebraic flux correction (AFC) methodology, which was introduced by Kuzmin and Turek in \cite{Kuzmin2002525} and refined in a series of publications \cite{Kuzmin2004131,KBOI,Kuzmin2006513,KuzminArticle,Kuzmin2010,KBNI,KBNII}. The family of AFC schemes is designed with the overall goal to prevent the creation of spurious oscillations by modifying the system matrix stemming from a standard Galerkin method in mass-conservative fashion. This algebraic design principle makes them particularly attractive for use in high-order isogeometric methods.

\section{High-resolution isogeometric analysis}
\label{sec:highres_iga}

This section briefly describes the basic construction principles of high-resolution isogeometric schemes for convection-dominated problems based on an extension of the AFC paradigm to B-Spline based discretizations of higher order. 

\subsection{Model problem}
\label{sec:high_res_iga:model_problem}

Consider the stationary convection-diffusion problem \cite{IGABook}
\begin{align}
- d \Delta u(\mathbf{x})) + \nabla \cdot (\mathbf{v} u(\mathbf{x})) &= 0 && \text{in} \; \Omega
\label{eq:model1}
\\
u(\mathbf{x}) &= \beta (\mathbf{x}) && \text{on} \; \varGamma
\label{eq:model2}
\end{align}
with diffusion coefficient $d = 0.0001$ and constant velocity vector $\mathbf{v}=[\sqrt{2}, \sqrt{2}]^\top$. The problem is solved on the two domains depicted in Figure~\ref{cnet}. 

Starting from the open knot vector $\Xi = [ 0,\, 0, \, 0, \, 0.5,\, 1,\, 1,\, 1]$, quadratic B-Spline basis functions $N_{a,2}(\xi)$ are generated by the Cox-de-Boor recursion formula \cite{deBoor1971}:
\begin{eqnarray}
    p=0:\quad N_{a,0}(\xi)&=&
    \left\{
    \begin{array}{ll}
    1 & \text{if } \xi_a\le\xi<\xi_{a+1},\\
    0 & \text{otherwise},
    \end{array}
    \right.\\
    p>0:\quad N_{a,p}(\xi)&=&
    \frac{\xi-\xi_a}{\xi_{a+p}-\xi_a}
    N_{a,p-1}(\xi)
    +
    \frac{\xi_{a+p+1}-\xi}{\xi_{a+p+1}-\xi_{a+1}}
    N_{a+1,p-1}(\xi),
\end{eqnarray}
where $\xi_a$ are the entries in the knot vector $\Xi$. Their tensor product construction yields the bivariate B-Spline basis functions $\hat\varphi_j(\xi,\eta)=N_a(\xi)N_b(\eta)$ (with index map $j\mapsto (a,b)$), which are used to define the computational geometry model
\begin{equation}
    \mathbf{x}(\xi,\eta)=\sum_j\mathbf{c}_j\hat\varphi_j(\xi,\eta),
    \quad
    (\xi,\eta)\in\hat\Omega=[0,1]^2
\end{equation}
with control points $\mathbf{c}_{j}\in\mathbb{R}^2$ indicated by dots in Fig.~\ref{cnet}. The mapping $\phi:\hat\Omega\to\Omega$ converts parametric values $\boldsymbol{\xi}=(\xi,\eta)$ into physical coordinates $\mathbf{x}=(x,y)$. The mapping should be bijective in order to possess a valid 'pull-back' operator $\phi^{-1}:\Omega\to\hat\Omega$.

For simplicity the boundary conditions are prescribed in the parametric domain:
\begin{equation}
\beta(\mathbf{x}=\phi(\boldsymbol{\xi})) = 
\begin{cases}
1  & \quad \text{if } \eta \leq \frac{1}{5}-\frac{1}{5}\xi\\
0  & \quad \text{otherwise}.
\end{cases}
\label{end1}
\end{equation}

\definecolor{rvwvcq}{rgb}{0,0,0}
\definecolor{cqcqcq}{rgb}{0.7529411764705882,0.7529411764705882,0.7529411764705882}
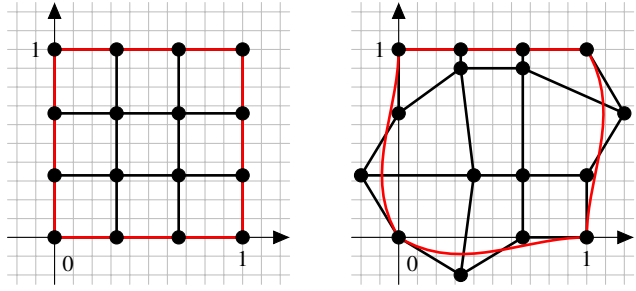
\begin{figure}[t!]
\centering
\begin{tikzpicture}[line cap=round,line join=round,>=triangle 45,x=2.5cm,y=2.5cm]
\draw [color=cqcqcq,, xstep=0.25cm,ystep=0.25cm] (-0.25,-0.25) grid (1.25,1.25);
\draw[->,color=black] (-0.25,0) -- (1.25,0);
\foreach \x in {1}
\draw[shift={(\x,0)},color=black] (0pt,2pt) -- (0pt,-2pt) node[below] {\footnotesize $\x$};
\draw[->,color=black] (0,-0.25) -- (0,1.25);
\foreach \y in {1}
\draw[shift={(0,\y)},color=black] (2pt,0pt) -- (-2pt,0pt) node[left] {\footnotesize $\y$};
\draw[color=black] (0pt,-10pt) node[right] {\footnotesize $0$};
\clip(-0.4,-0.4) rectangle (1.4,1.4);
\draw [line width=1pt] (0,0.33)-- (0,0.66);
\draw [line width=1pt] (0,0.66)-- (0,1);
\draw [line width=1pt] (0,1)-- (0.33,1);
\draw [line width=1pt] (0.33,1)-- (0.66,1);
\draw [line width=1pt] (0.66,1)-- (1,1);
\draw [line width=1pt] (1,1)-- (1,0.66);
\draw [line width=1pt] (1,0.66)-- (1,0.33);
\draw [line width=1pt] (1,0.33)-- (1,0);
\draw [line width=1pt] (1,0)-- (0.66,0);
\draw [line width=1pt] (0.66,0)-- (0.33,0);
\draw [line width=1pt] (0.33,0)-- (0,0);
\draw [line width=1pt] (0,0)-- (0,0.33);
\draw [line width=1pt] (0,0.33)-- (0.33,0.33);
\draw [line width=1pt] (0.33,0.33)-- (0.66,0.33);
\draw [line width=1pt] (0.66,0.33)-- (1,0.33);
\draw [line width=1pt] (0,0.66)-- (0.33,0.66);
\draw [line width=1pt] (0.33,0.66)-- (0.66,0.66);
\draw [line width=1pt] (0.66,0.66)-- (1,0.66);
\draw [line width=1pt] (0.66,1)-- (0.66,0.66);
\draw [line width=1pt] (0.66,0.66)-- (0.66,0.33);
\draw [line width=1pt] (0.66,0.33)-- (0.66,0);
\draw [line width=1pt] (0.33,0)-- (0.33,0.33);
\draw [line width=1pt] (0.33,0.33)-- (0.33,0.66);
\draw [line width=1pt] (0.33,0.66)-- (0.33,1);

\draw [red, line width=1pt] (0,0)-- (1,0);
\draw [red, line width=1pt] (0,0)-- (0,1);
\draw [red, line width=1pt] (0,1)-- (1,1);
\draw [red, line width=1pt] (1,0)-- (1,1);

\begin{scriptsize}
\draw [fill=rvwvcq] (0,0) circle (2.5pt);
\draw [fill=rvwvcq] (0.33,0) circle (2.5pt);
\draw [fill=rvwvcq] (0.66,0) circle (2.5pt);
\draw [fill=rvwvcq] (1,0) circle (2.5pt);
\draw [fill=rvwvcq] (0,0.33) circle (2.5pt);
\draw [fill=rvwvcq] (0.33,0.33) circle (2.5pt);
\draw [fill=rvwvcq] (0.66,0.33) circle (2.5pt);
\draw [fill=rvwvcq] (1,0.33) circle (2.5pt);
\draw [fill=rvwvcq] (0,0.66) circle (2.5pt);
\draw [fill=rvwvcq] (0.33,0.66) circle (2.5pt);
\draw [fill=rvwvcq] (0.66,0.66) circle (2.5pt);
\draw [fill=rvwvcq] (1,0.66) circle (2.5pt);
\draw [fill=rvwvcq] (0,1) circle (2.5pt);
\draw [fill=rvwvcq] (0.33,1) circle (2.5pt);
\draw [fill=rvwvcq] (0.66,1) circle (2.5pt);
\draw [fill=rvwvcq] (1,1) circle (2.5pt);
\end{scriptsize}
\end{tikzpicture}	
\begin{tikzpicture}[line cap=round,line join=round,>=triangle 45,x=2.5cm,y=2.5cm]

\draw [color=cqcqcq,, xstep=0.25cm,ystep=0.25cm] (-0.25,-0.25) grid (1.25,1.25);
\draw[->,color=black] (-0.25,0) -- (1.25,0);
\foreach \x in {1}
\draw[shift={(\x,0)},color=black] (0pt,2pt) -- (0pt,-2pt) node[below] {\footnotesize $\x$};
\draw[->,color=black] (0,-0.25) -- (0,1.25);
\foreach \y in {1}
\draw[shift={(0,\y)},color=black] (2pt,0pt) -- (-2pt,0pt) node[left] {\footnotesize $\y$};
\draw[color=black] (0pt,-10pt) node[right] {\footnotesize $0$};
\clip(-0.4,-0.4) rectangle (1.4,1.4);
\draw [line width=1pt] (-0.2,0.33)-- (0,0.66);
\draw [line width=1pt] (0,0.66)-- (0,1);
\draw [line width=1pt] (0,1)-- (0.33,1);
\draw [line width=1pt] (0.33,1)-- (0.66,1);
\draw [line width=1pt] (0.66,1)-- (1,1);
\draw [line width=1pt] (1,1)-- (1.2,0.66);
\draw [line width=1pt] (1.2,0.66)-- (1,0.33);
\draw [line width=1pt] (1,0.33)-- (1,0);
\draw [line width=1pt] (1,0)-- (0.66,0);
\draw [line width=1pt] (0.66,0)-- (0.33,-0.2);
\draw [line width=1pt] (0.33,-0.2)-- (0,0);
\draw [line width=1pt] (0,0)-- (-0.2,0.33);
\draw [line width=1pt] (-0.2,0.33)-- (0.4,0.33);
\draw [line width=1pt] (0.4,0.33)-- (0.66,0.33);
\draw [line width=1pt] (0.66,0.33)-- (1,0.33);
\draw [line width=1pt] (0,0.66)-- (0.33,0.9);
\draw [line width=1pt] (0.33,0.9)-- (0.66,0.9);
\draw [line width=1pt] (0.66,0.9)-- (1.2,0.66);
\draw [line width=1pt] (0.66,1)-- (0.66,0.9);
\draw [line width=1pt] (0.66,0.9)-- (0.66,0.33);
\draw [line width=1pt] (0.66,0.33)-- (0.66,0);
\draw [line width=1pt] (0.33,-0.2)-- (0.4,0.33);
\draw [line width=1pt] (0.4,0.33)-- (0.33,0.9);
\draw [line width=1pt] (0.33,0.9)-- (0.33,1);

\draw [red, line width=1pt] plot (0,0) .. controls (-0.2,0.33) and (0,0.66) .. (0,1);
\draw [red, line width=1pt] plot (0,0) .. controls (0.33,-0.2) and (0.66,0) .. (1,0);
\draw [red, line width=1pt] plot (1,0) .. controls (1.0,0.33) and (1.2,0.66) .. (1,1);
\draw [red, line width=1pt] plot (0,1) .. controls (0.33,1.00) and (0.66,1.0) .. (1,1);

\begin{scriptsize}
\draw [fill=rvwvcq] (0,0) circle (2.5pt);
\draw [fill=rvwvcq] (0.33,-0.2) circle (2.5pt);
\draw [fill=rvwvcq] (0.66,0) circle (2.5pt);
\draw [fill=rvwvcq] (1,0) circle (2.5pt);
\draw [fill=rvwvcq] (-0.2,0.33) circle (2.5pt);
\draw [fill=rvwvcq] (0.4,0.33) circle (2.5pt);
\draw [fill=rvwvcq] (0.66,0.33) circle (2.5pt);
\draw [fill=rvwvcq] (1,0.33) circle (2.5pt);
\draw [fill=rvwvcq] (0,0.66) circle (2.5pt);
\draw [fill=rvwvcq] (0.33,0.9) circle (2.5pt);
\draw [fill=rvwvcq] (0.66,0.9) circle (2.5pt);
\draw [fill=rvwvcq] (1.2,0.66) circle (2.5pt);
\draw [fill=rvwvcq] (0,1) circle (2.5pt);
\draw [fill=rvwvcq] (0.33,1) circle (2.5pt);
\draw [fill=rvwvcq] (0.66,1) circle (2.5pt);
\draw [fill=rvwvcq] (1,1) circle (2.5pt);
\end{scriptsize}
\end{tikzpicture}

\caption{Unit square (left) and deformed domain (right) modeled by tensor-product quadratic B-Spline basis functions defined on the open knot vector $\Xi = [ 0,\, 0, \, 0, \, 0.5,\, 1,\, 1,\, 1]$.}
\label{cnet}
\end{figure}

\subsection{Galerkin method} 
\label{sec:highres_iga:galerkin}

Application of the Galerkin method to \eqref{eq:model1}--\eqref{end1} yields: Find $u^h\in S^h$ such that
\begin{equation}
	d \int_\Omega \nabla u^h \cdot \nabla v^h \mathrm{d}\mathbf{x}
	+ 
	\int_\Omega \nabla\cdot  (\mathbf{{v}}u)^h v^h) \mathrm{d}\mathbf{x} 
	=
	\int_\Omega R v^h \mathrm{d}\mathbf{x}
	\label{eq:galerkin}
\end{equation}
for all test functions $v^h\in V^h$ that vanish on the entire boundary $\Gamma$ due to the prescription of Dirichlet boundary conditions. In the framework of IGA the  discrete spaces $S^h$ and $V^h$ are spanned by multivariate B-Spline basis functions $\lbrace\varphi_j(\mathbf{x})\rbrace$. 

Using Fletcher's group formulation \cite{Fletcher}, the approximate solution $u^h$ and the convective flux $(\mathbf{{v}}u)^h$ can be represented as follows \cite{KBNI}: 
\begin{equation}
	u^h(\mathbf{x}) = \sum_j u_j \varphi_j(\mathbf{x}),\quad
	(\mathbf{{v}}u)^h(\mathbf{x}) = \sum_j (\mathbf{{v}}_j u_j) \varphi_j(\mathbf{x}).
	\label{eq:fletcher}
\end{equation}
Substitution into \eqref{eq:galerkin} and replacing $v^h$ by all basis functions yields the matrix form 
\begin{equation}
	(S-K) \mathbf{u} = \mathbf{r},
	\label{scheme}
\end{equation}
where $\mathbf{u}$ is the vector of coefficients $u_i$ used in the expansion of the solution \eqref{eq:fletcher} and the entries of the discrete diffusion ($S = \lbrace s_{ij}\rbrace$) and convection ($K = \lbrace k_{ij}\rbrace$) operators and the discretized right-hand side vector ($\mathbf{r} = \lbrace r_{i}\rbrace$) are given by
\begin{align}
    k_{ij} &= -\mathbf{v}_j \cdot \mathbf{c}_{ij},
	&
	\mathbf{c}_{ij} &=\int_\Omega \nabla \varphi_j \varphi_i \mathrm{d}\mathbf{x},
	\\
	s_{ij} &= d\int_\Omega \nabla \varphi_j \cdot \nabla \varphi_i\mathrm{d}\mathbf{x},
	&
	r_i &= \int_\Omega R \varphi_i \mathrm{d}\mathbf{x}.
\end{align}

The above integrals are assembled by resorting to numerical quadrature over the unit square $\hat\Omega=[0,1]^2$ using the 'pull-back' operator $\phi^{-1}:\Omega\to\hat\Omega$. For the entries of the physical diffusion matrix the final expression reads as follows \cite{Jaeschke}:
\begin{equation}
    s_{ij}
    =
    d\int_{\hat\Omega}
    \nabla_{\boldsymbol{\xi}}\hat\varphi_j(\boldsymbol{\xi})\cdot
    G(\boldsymbol{\xi})
    \nabla_{\boldsymbol{\xi}}\hat\varphi_j(\boldsymbol{\xi})
    \mathrm{d}\boldsymbol{\xi},
    \label{eq:anisodiff}
\end{equation}
where the geometric factor $G(\boldsymbol{\xi})$ is given in terms of the Jacobian $J=D\phi$:
\begin{eqnarray}
    G(\boldsymbol{\xi})
    =
    |\det J(\boldsymbol{\xi})|
    J^{-1}(\boldsymbol{\xi})
    J^{-\top}(\boldsymbol{\xi}).
    \label{eq:geofac}
\end{eqnarray}
It should be noted that expression \eqref{eq:anisodiff} can be interpreted as the discrete counterpart of an anisotropic diffusion problem with symmetric diffusion tensor $dG(\boldsymbol{\xi})$ that is solved on the unit square $\hat\Omega$ using tensor-product B-Splines on a perpendicular grid.

\subsection{Algebraic flux correction}
\label{sec:highres_iga:afc}

The isogeometric Galerkin method \eqref{scheme} is turned into a stabilized high-resolution scheme by applying the principles of algebraic flux correction (AFC) of TVD-type, which were developed for lowest-order Lagrange finite elements in \cite{KuzminArticle,Kuzmin2004131}.

In essence, the discrete convection operator $K$ is modified in two steps:
\begin{itemize}
    \item[1.] Eliminate negative off-diagonal entries from  $K$ by adding a \emph{discrete diffusion operator} $D$ to obtain the modified discrete convection operator $L=K+D$.
    \item[2.] Remove excess artificial diffusion in regions where this is possible without generating spurious wiggles by applying \emph{non-linear anti-diffusion}: $K^*(\mathbf{u})=L+\bar{F}(\mathbf{u})$.
\end{itemize}
\textbf{Discrete diffusion operator.} The optimal entries of $D = \lbrace d_{ij} \rbrace$ are given by \cite{Kuzmin2002525}:
\begin{equation}
d_{ij} = d_{ji} = \max \lbrace 0, -k_{ij}, -k_{ji} \rbrace,\quad
d_{ii} = - \sum_{j \neq i} d_{ij},
\label{dij}
\end{equation}
yielding a symmetric operator with zero column and row sums. The latter enables the decomposition of the diffusive contribution to the $i$th degree of freedom
\begin{equation}
\left(D\mathbf{u}\right)_i = \sum_{j \neq i} f_{ij},
\quad
f_{ij} = d_{ij}(u_j-u_i),
\label{dec}
\end{equation}
whereby the diffusive fluxes $f_{ij}=-f_{ji}$ are skew-symmetric by design \cite{Kuzmin2002525}. 

In a practical implementation, operator $D$ is not constructed explicitly, but the entries of $L:=K$ are modified in a loop over all pairs of degrees of freedoms $(i,j)$ for which $j\ne i$ and the basis functions have overlapping support $\text{supp}\hat\varphi_i\cap\text{supp}\hat\varphi_j\ne \emptyset$. For univariate B-Spline basis functions of order $p$, we have $\text{supp}\hat\varphi_i=(\xi_i,\xi_{i+p+1})$, where $\xi_i$ denotes the $i$th entry of the knot vector $\Xi$. Hence, the loops in \eqref{dij} and \eqref{dec} extend over all $j\ne i$ with $|j-i|\le p$ in one spatial dimension, which can be easily generalized to tensor-product B-Splines in multiple dimensions.

The modified convection operator $L=K+D$ yields the  stabilized linear scheme
\begin{equation}
(S-L) \mathbf{u} = \mathbf{r}.
\label{art}
\end{equation}

\textbf{Anisotropic physical diffusion.} The discrete diffusion matrix $S$ might also cause spurious oscillations in the solution since it is only 'harmless' for lowest order finite elements under the additional constraint that triangles are nonobtuse (all angles smaller than or equal to $\pi/2$) and quadrilaterals are nonnarrow (aspect ratios smaller than or equal to $\sqrt{2}$ \cite{Farago}), respectively. Kuzmin et al. \cite{Kuzmin2008,KBNI} propose stabilization techniques for anisotropic diffusion problems, which can be applied to \eqref{eq:anisodiff} directly. It should be noted, however, that we did not observe any spurious wiggles in all our numerical tests even without any special treatment of the diffusion matrix $S$.

\textbf{Nonlinear anti-diffusion.} According to Godunov's theorem \cite{Godunov}, the linear scheme \eqref{art} is limited to first-order accuracy. Therefore a nonlinear scheme must be constructed by adaptively blending between schemes (\ref{scheme}) and (\ref{art}), namely \cite{KuzminArticle,Kuzmin2004131}:
\begin{equation}
(S-K^*(\mathbf{u})) \mathbf{u} = \mathbf{r}.
\label{sch}
\end{equation}
Here, the nonlinear discrete convection operator reads
\begin{equation} 
K^*(\mathbf{u}) = L + \bar{F}(\mathbf{u}) = K + D + \bar{F}(\mathbf{u}),
\end{equation}
which amounts to applying a modulated anti-diffusion operator $\bar{F}(\mathbf{u})$ to avoid the loss of accuracy in smooth regions due to excessive artificial diffusion. The raw anti-diffusion, $-D$, features all properties of a discrete diffusion operator, and hence, its contribution to a single degree of freedom can be decomposed as follows \cite{KuzminArticle,Kuzmin2004131}:
\begin{equation}
f_i(\mathbf{u})
:=
(\bar{F}(\mathbf{u})\mathbf{u})_i
=
\sum_{j \neq i} \alpha_{ij}(\mathbf{u})d_{ij}(u_i-u_j),
\end{equation}
where $\alpha_{ij}(\mathbf{u})=\alpha_{ji}(\mathbf{u})$ is an adaptive flux limiter. Clearly, for $\alpha_{ij}\equiv 1$ the anti-diffusive fluxes will restore the original Galerkin scheme (\ref{scheme}) and $\alpha_{ij}\equiv 0$ will lead to the linear scheme \eqref{art}. Kuzmin et al. \cite{KuzminArticle,Kuzmin2004131} proposed a TVD-type multi-dimensional limiting strategy for lowest-order Lagrange finite elements, which ensures that the resulting scheme \eqref{sch} yields accurate solutions that are free of spurious oscillations. The flux limiter was extended to non-nodal basis functions in \cite{Jaeschke} and utilized for computing the numerical results presented in Section~\ref{sec:results}.


Like with the diffusion operator $D$, we do not construct $K^*(\mathbf{u})$ explicitly but include the anti-diffusive correction $\mathbf{\bar f}(\mathbf{u})=\lbrace \bar f_i(\mathbf{u}) \rbrace$ into the right-hand side \cite{KuzminArticle,Kuzmin2004131}
\begin{equation}
(S-L) \mathbf{u} = \mathbf{r} + \mathbf{\bar f}(\mathbf{u}).
\label{sch2}
\end{equation}
The nonlinear scheme can be solved by iterative defect correction \cite{Kuzmin2004131} possibly combined with Anderson acceleration \cite{Walker2011} or by an inexact Newton method \cite{Moller2007}.

Non-nodal degrees of freedom make it necessary to first project the prescribed boundary values \eqref{end1} onto the solution space $S^h$ so that the coefficients of the degrees of freedoms that are located at the Dirichlet boundary part can be overwritten accordingly. Since the standard $L_2$ projection can lead to non-physical under- and overshoots near discontinuities and steep gradients, the constrained data projection approach proposed in \cite{Kuzmin2010} for lowest order nodal finite elements is used.

\section{Numerical results}
\label{sec:results}

This section presents the numerical results for the model problem \eqref{eq:model1}--\eqref{end1}, which were computed using the open-source isogeometric analysis library G+Smo \cite{jlmmz2014}. 

The tensor-product B-Spline basis ($4\times 4$ basis functions of degree $p=2$) that was used for the geometry models depicted in Fig.~\ref{cnet} was refined by means of knot insertion \cite{IGABook} to generate $18 \times 18$ quadratic B-Spline basis functions for approximating the solution. It should be noted that this type of refinement, which is an integral part of the Isogeometric Analysis framework, preserves the shape of the geometry exactly. Consequently, the numerical solution does not suffer from an additional error stemming from an approximated computational domain as it is the case for, say, higher-order Lagrange finite elements defined on simplex or quadrilateral meshes.

\begin{figure}[b!]
\centering
\includegraphics[height=4.2cm, trim=150 0 0 60, clip]{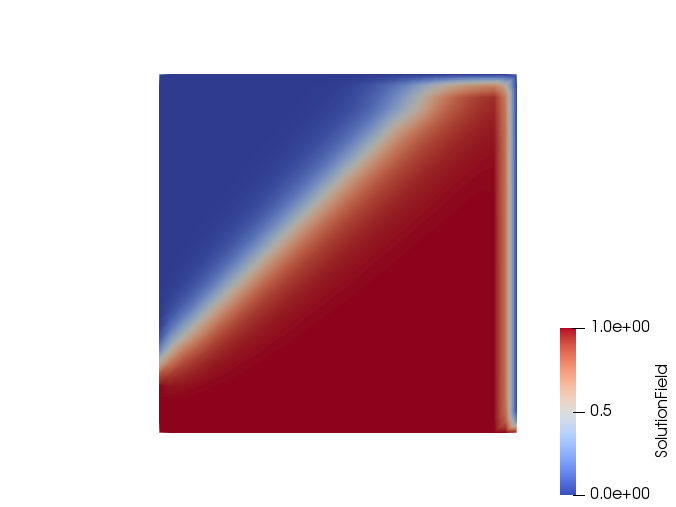}\quad
\includegraphics[height=4.2cm, trim=100 0 0 60, clip]{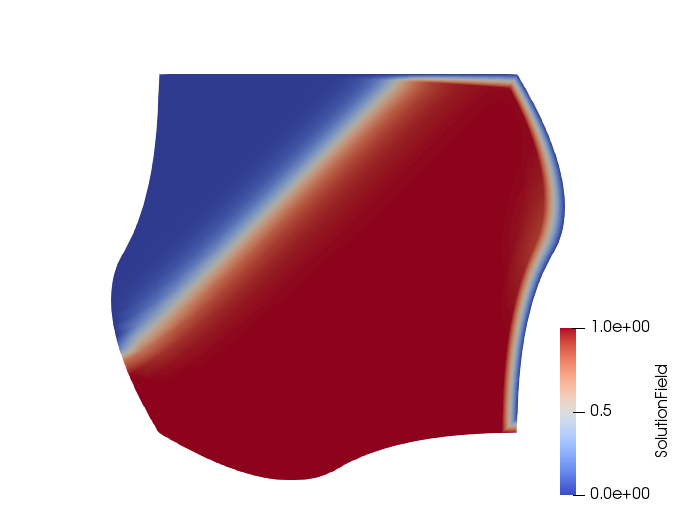}
\caption{Numerical solutions computed on the unit square (left) and deformed domain (right).}
\label{original_result}
\end{figure}

For the diffusion coefficient $d=0.0001$ and the considered basis the element P{\'e}clet number is equal to $Pe_h \approx 555$, which states that the problem is highly convection-dominated. For the deformed geometry, the actual value varies slightly from one 'element' to the other but stays in the same order of magnitude.


The numerical solution that was computed on the unit square is depicted in Fig. \ref{original_result} (left), whereas the approximate solution for the deformed geometry is shown on the right. In both cases the minimum and maximum bounds of the exact solution, that is, $u_{min}=0$ and $u_{max}=1$ are preserved by the numerical counterpart, which results from the successful application of the AFC stabilization of TVD-type.

It should be noted that the internal layer is smeared stronger than the boundary layer, which is due to the constrained $L_2$ projection of the Dirichlet boundary data into the space $S^h$. The discontinuous profile \eqref{end1} along the left boundary cannot be represented exactly by quadratic B-Splines, and hence, it is smeared across multiple 'elements'. A possible remedy is to locally reduce the approximation order to $p=1$ by inserting a knot at the boundary location $\eta_b=1/5$ and increasing its multiplicity to $m_b=2$, which will reduce the continuity to $C^{p-m_{b}}=C^0$ locally. The varying thickness of the boundary layer on the deformed geometry stems from the fact the distance of the rightmost vertical internal 'grid line' to the boundary also varies.


\section{Conclusions}
\label{sec:conclusions}

The high-resolution isogeometric scheme presented in this work for the stationary convection-diffusion equation is a first step to establish isogeometric methods for convection-dominated problems and, in particular, compressible flows, which are addressed in more detail in \cite{igaAFC2_fef2017}. This chapter extends the family of algebraic flux correction schemes to quadratic B-Spline discretizations thereby demonstrating that the algebraic design principles that were originally derived for low-order nodal Lagrange finite elements carry over to non-nodal Spline basis functions.

Ongoing research focuses on the extension of this approach to truncated hierarchical B-Splines \cite{Giannelli2012} possibly combined with the local increase of the knot multiplicities, which seems to be a viable approach for refining the spline spaces $S^h$ and $V^h$ adaptively in the vicinity of shocks and steep gradients to compensate for the local reduction of the approximation order by algebraic flux correction ($h$-refinement) and to prevent excessive spreading of these localized features (continuity reduction).

\begin{acknowledgement}
This work has been supported by the European Unions Horizon 2020 research and innovation programme under grant agreement No. 678727.
\end{acknowledgement}

\bibliographystyle{unsrt}
\bibliography{bib}

\begin{thebibliography}{10}

\bibitem{igaAFC2_fef2017}
M.~M\"oller and A.~Jaeschke.
\newblock High-order isogeometric methods for compressible flows. {II}.
  {C}ompressible {E}uler equations, Accepted for publication, online version
  available at arXiv.

\bibitem{IGAArticle}
T.J.R. Hughes, J.A. Cottrell, and Y.~Bazilevs.
\newblock Isogeometric analysis: {CAD}, finite elements, {NURBS}, exact
  geometry and mesh refinement.
\newblock {\em Computer Methods in Applied Mechanics and Engineering},
  194(39–41):4135 -- 4195, 2005.

\bibitem{IGABook}
J.A. Cottrell, T.J.R. Hughes, and Y.~Bazilevs.
\newblock {\em Isogeometric Analysis: Toward Integration of {CAD} and {FEA}}.
\newblock John Wiley \& Sons, Ltd., 2009.

\bibitem{IGAart2}
P.~Trontin.
\newblock Isogeometric analysis of {E}uler compressible flow. {A}pplication to
  aerodynamics.
\newblock {\em 50th AIAA}, 2012.

\bibitem{Jaeschke}
A.~Jaeschke.
\newblock Isogeometric analysis for compressible flows with application in
  turbomachinery.
\newblock Master's thesis, Delft University of Technology, 2015.

\bibitem{Brooks1982199}
A.N. Brooks and T.J.R. Hughes.
\newblock Streamline upwind/{P}etrov-{G}alerkin formulations for convection
  dominated flows with particular emphasis on the incompressible
  {N}avier-{S}tokes equations.
\newblock {\em Computer Methods in Applied Mechanics and Engineering},
  32(1–3):199 -- 259, 1982.

\bibitem{Kuzmin2002525}
D.~Kuzmin and S.~Turek.
\newblock Flux correction tools for finite elements.
\newblock {\em Journal of Computational Physics}, 175(2):525 -- 558, 2002.

\bibitem{Kuzmin2004131}
D.~Kuzmin and S.~Turek.
\newblock High-resolution {FEM}-{TVD} schemes based on a fully multidimensional
  flux limiter.
\newblock {\em Journal of Computational Physics}, 198(1):131 -- 158, 2004.

\bibitem{KBOI}
D.~Kuzmin and M.~M{\"o}ller.
\newblock {\em Flux-Corrected Transport}, chapter Algebraic flux correction I.
  {S}calar conservation laws.
\newblock Springer, 2005.

\bibitem{Kuzmin2006513}
D.~Kuzmin.
\newblock On the design of general-purpose flux limiters for finite element
  schemes. {I}. {S}calar convection.
\newblock {\em Journal of Computational Physics}, 219(2):513 -- 531, 2006.

\bibitem{KuzminArticle}
D.~Kuzmin.
\newblock Algebraic flux correction for finite element discretizations of
  coupled systems.
\newblock {\em Computational Methods for Coupled Problems in Science and
  Engineering II, CIMNE, Barcelona}, pages 653--656, 2007.

\bibitem{Kuzmin2010}
D.~Kuzmin, M.~M{\"o}ller, J.N. Shadid, and M.~Shashkov.
\newblock Failsafe flux limiting and constrained data projections for equations
  of gas dynamics.
\newblock {\em Journal of Computational physics}, 229(23):8766--8779, 11 2010.

\bibitem{KBNI}
D.~Kuzmin.
\newblock {\em Flux-Corrected Transport}, chapter Algebraic flux correction I.
  Scalar conservation laws.
\newblock Springer, 2012.

\bibitem{KBNII}
D.~Kuzmin, M.~M{\"o}ller, and M.~Gurris.
\newblock {\em Flux-Corrected Transport}, chapter Algebraic flux correction II.
  Compressible flow problems.
\newblock Springer, 2012.

\bibitem{deBoor1971}
C.~de~Boor.
\newblock Subroutine package for calculating with {B}-splines.
\newblock Technical Report LA-4728-MS, Los Alamos Scient. Lab., 1971.

\bibitem{Fletcher}
C.A.J. Fletcher.
\newblock The group finite element formulation.
\newblock {\em Comput. Methods Appl. Mech. Eng}, 37:225--243, 1983.

\bibitem{Farago}
I.~Farago, R.~Horvath, and S.~Korotov.
\newblock Discrete maximum principle for linear parabolic problems solved on
  hybrid meshes.
\newblock {\em Applied Numerical Mathematics}, 53(2-4):249--264, 05 2004.

\bibitem{Kuzmin2008}
D.~Kuzmin, M.J. Shashkov, and D.~Svyatski.
\newblock A constrained finite element method satisfying the discrete maximum
  principle for anisotropic diffusion problems on arbitrary meshes, 2008.

\bibitem{Godunov}
S.K. Godunov.
\newblock Finite difference method for numerical computation of discontinous
  solutions of the equations of fluid dynamics.
\newblock {\em Mat. Sbornik}, 47:271--306, 1959.

\bibitem{Walker2011}
H.F. Walker and P.~Ni.
\newblock Anderson acceleration for fixed-point iterations.
\newblock {\em SIAM J. Numer. Anal.}, 49(4):1715--1735, 08 2011.

\bibitem{Moller2007}
M.~M{\"o}ller.
\newblock Efficient solution techniques for implicit finite element schemes
  with flux limiters.
\newblock {\em International Journal for Numerical Methods in Fluids},
  55(7):611--635, 11 2007.

\bibitem{jlmmz2014}
B.~J{\"u}ttler, U.~Langer, A.~Mantzaflaris, S.~Moore, and W.~Zulehner.
\newblock Geometry + {S}imulation {M}odules: Implementing isogeometric
  analysis.
\newblock {\em Proc. Appl. Math. Mech.}, 14(1):961--962, 2014.
\newblock Special Issue: 85th Annual Meeting of the Int. Assoc. of Appl. Math.
  and Mech. (GAMM), Erlangen 2014.

\bibitem{Giannelli2012}
C.~Giannelli, B.~J\"{u}ttler, and H.~Speleers.
\newblock {THB}-splines: The truncated basis for hierarchical splines.
\newblock {\em Computer Aided Geometry Design}, 29(7):485--498, 10 2012.

\end{thebibliography}
\end{document}